\documentclass[a4paper,11pt]{smfart}

\usepackage{mathpazo}

\usepackage{booktabs}
\usepackage{amsmath, amsfonts, amssymb, amsthm}
\usepackage[all]{xy}
\usepackage[height=22.5cm, width=15.5cm]{geometry}

\providecommand{\og}{``}
\providecommand{\fg}{''}
\providecommand{\smfandname}{\&}

\numberwithin{equation}{section}

\theoremstyle{plain}
\newtheorem{thm}{Th\'eor\`eme}[section]
\newtheorem{prop}[thm]{Proposition}
\newtheorem{lem}[thm]{Lemme}

\newtheorem*{thm*}{Th\'eor\`eme}

\theoremstyle{definition}
\newtheorem{defn}[thm]{D\'efinition}

\theoremstyle{remark}
\newtheorem*{rem}{Remarque}
\newtheorem*{ex}{Exemple}

\newcommand{\mbb}[1]{\mathbb{#1}}
\newcommand{\lie}[1]{{\mathfrak{#1}}}
\newcommand{\abs}[1]{\lvert #1\rvert}

\DeclareMathOperator{\Aut}{Aut}

\title{Sur les quotients discrets de semi-groupes complexes}
\author{Christian Miebach}\thanks{L'auteur est soutenu par SFB/TR~12 du DFG}
\address{Fakult\"at f\"ur Mathematik, Ruhr-Universit\"at Bochum,
Universit\"atsstra{\ss}e 150, D - 44780 Bochum}
\email{christian.miebach@ruhr-uni-bochum.de}
\address{Adresse actuelle: Centre de Math\'ematiques et d'Informatique,
UMR-CNRS 6632 (LATP), 39, rue Joliot-Curie, Universit\'e de Provence,
13453 Marseille Cedex 13 France}
\date{22 Avril 2008}

\begin{document}

\maketitle

\begin{abstract}
Soit $X=G/K$ un espace sym\'etrique hermitien irr\'educible de type non-compact
et soit $S\in G^\mbb{C}$ le semi-groupe associ\'e form\'e des compressions de
$X$. Soit $\Gamma\subset G$ un sous-groupe discret. Nous donnons une condition
suffisante pour que le quotient $\Gamma\backslash S$ soit une vari\'et\'e
de Stein. En outre nous d\'emontrons qu'en g\'en\'eral $\Gamma\backslash S$
n'est pas de Stein ce qui r\'efute une conjecture de Achab, Betten et Kr\"otz.
\end{abstract}

\section*{Introduction}

Soit $X=G/K$ un espace sym\'etrique hermitien irr\'educible de type non-compact
et soit $S\in G^\mbb{C}$ le semi-groupe associ{\'e} consistant des compressions
de $X$. On d\'emontre dans~\cite{AchBeKr} que le quotient $\Gamma\backslash S$
est holomorphiquement s\'eparable pour tout sous-groupe discret $\Gamma\subset
G$. De plus, ils conjecturent que ce quotient est une vari\'et\'e de Stein et
ils v\'erifient cette conjecture dans le cas o\`u $G={\rm{SL}}(2,\mbb{R})$ et
$\Gamma\subset{\rm{SL}}(2,\mbb{Z})$.

Dans ce travail nous donnons une condition suffisante pour que $\Gamma\backslash
S$ soit de Stein.

\begin{thm*}
Soit $X=G/K$ un espace sym\'etrique Hermitien irr\'educible de type non-compact
et soit $S\in G^\mbb{C}$ le semi-groupe associ{\'e} form\'e des compressions
de $X$. Soit $\Gamma\subset G$ un sous-groupe discret qui agit librement sur
$X$. Si $X/\Gamma$ est une vari\'et\'e de Stein, alors $\Gamma\backslash S$ est 
de Stein aussi.
\end{thm*}

Apr\`es avoir pr\'esent\'e les fondements necessaires de la th\'eorie des
espaces sym\'etriques hermitiens et des semi-groupes de compression dans la
premi\`ere section, nous d\'emontrons ce th\'eor\`eme. Notre preuve utilise
essentiellement le fibr\'e $G^\mbb{C}\to G^\mbb{C}/K^\mbb{C}\to Z$ (o\`u nous
notons par $Z$ le dual compact de $X$) et des techniques qui ont \'et\'e
d\'evelopp\'ees dans~\cite{Mie4}. Dans la derni\`ere section nous montrons
que, si $G={\rm{SL}}(2,\mbb{R})$ et si $G/\Gamma$ est compacte, alors
$\Gamma\backslash S$ n'est pas de Stein ce qui r\'efute la conjecture
g\'en\'erale dans~\cite{AchBeKr}.

Je remercie Karl~Oeljeklaus pour l'int\'er\^et et l'aide qu'il a
apport\'e \`a ce travail ainsi que pour l'invitation \`a LATP de l'Universit\'e
de Provence o\`u ce travail a \'et\'e effectu\'e.

\section{Le semi-groupe de compression associ\'e \`a un espace sym\'etrique
Hermitien}

Dans cette section nous parcourons les parties de la th\'eorie des espaces
sym\'etriques hermitiens non-compacts $X$ dont nous avons besoin afin de
d\'efinir le semi-groupe de compression associ\'e \`a $X$. Pour en savoir plus
nous prions le lecteur de s'adresser \`a~\cite{Hel} et~\cite{HiNe}.

Soit $X=G/K$ un espace sym\'etrique hermitien irr\'educible de type non-compact.
Nous pouvons supposer que $G$ est un groupe de Lie r\'eel connexe simple sans
facteurs compacts qui est plong\'e dans sa complexification universelle
$G^\mbb{C}$ et que $K$ est un sous-groupe maximal compact de $G$ qui est
d\'efini par l'involution de Cartan $\theta\in\Aut(G)$. Nous d\'enotons par
$\lie{g}=\lie{k}\oplus\lie{p}$ la decomposition de l'alg\`ebre de Lie $\lie{g}$
de $G$ par rapport \`a $\theta_*\in\Aut(\lie{g})$. Puisque $X$ est une
vari\'et\'e complexe sur laquelle $G$ agit par des transformations holomorphes,
il y a une structure complexe $J$ sur $\lie{p}$ invariante par l'action adjointe
de $K$. Nous \'etendons $J$ comme un operateur $\mbb{C}$--lin\'eaire \`a
$\lie{p}^\mbb{C}$. Puisque les valeurs propres de $J$ sont $\pm i$, on obtient
la d\'ecomposation en somme direct $\lie{p}^\mbb{C}=\lie{p}^+\oplus\lie{p}^-$.
On peut d\'emontrer que $\lie{p}^\pm$ est une sous-alg\`ebre commutative de
$\lie{g}^\mbb{C}$ et $[\lie{k}^\mbb{C},\lie{p}^\pm]\subset\lie{p}^\pm$.

Soit $K^\mbb{C}$ et $P^{\pm}$ les sous-groupes analytiques de $G^\mbb{C}$ dont
l'alg\`ebre de Lie est donn\'ee par $\lie{k}^\mbb{C}$ et $\lie{p}^{\pm}$.
Alors $K^\mbb{C}$ et $P^\pm$ sont des sous-groupes ferm\'es et complexes de
$G^\mbb{C}$. De plus, l'ensemble $Q^\pm:=K^\mbb{C}P^\pm$ est un sous-groupe
ferm\'e complexe parabolique et isomorphe au produit semi-direct
$K^\mbb{C}\ltimes P^\pm$. L'espace complexe homog\`ene $Z:=G^\mbb{C}/Q^-$ peut
\^etre identifi\'e au dual compact de $X$ et l'application $X\to Z$,
$gK\mapsto gQ^-$, d\'efinit un plongement holomorphe de $X$ comme un ouvert de
$Z$. Ce plongement est connu comme le plongement de Borel. D\'esormais
nous regardons $X$ comme l'orbite ouverte $G\cdot eQ^-\subset Z$.

\begin{defn}
L'ensemble
\begin{equation*}
S:=\bigl\{g\in G^\mbb{C};\ g(\overline{X})\subset X\bigr\}
\end{equation*}
s'appelle le semi-groupe ouvert complexe de compression associ\'e \`a $X$.
\end{defn}

Il suit imm\'ediatement que $S$ est un sous-ensemble ouvert et $(G\times
G)$--invariant de $G^\mbb{C}$ o\`u l'action de $G\times G$ est donn\'ee par
$z\mapsto g_1zg_2^{-1}$. De plus, $S$ est stable par rapport \`a la
multiplication de $G^\mbb{C}$ et par cons\'equent un semi-groupe. Ce
semi-groupe de compression est un cas particulier d'un semi-groupe
d'Ol'shanski{\u\i} (comparer chapitre~8 dans~\cite{HiNe}).

\begin{thm}[\cite{Ne2}]
Le semi-groupe de compression $S$ est un domaine d'holomorphie hyperbolique au
sens de Kobayashi dans $G^\mbb{C}$. En particulier, le groupe $G\times G$ agit
proprement sur $S$.
\end{thm}

Le quotient $\Gamma\backslash S$ est une vari\'et\'e complexe pour tout
sous-groupe $\Gamma\subset G$ qui agit par multiplication de gauche sur $S$.
Le th\'eor\`eme suivant est un cas particulier d'un resultat plus g\'en\`eral
de Achab, Betten et Kr\"otz (\cite{AchBeKr}).

\begin{thm}\label{Thm:HolSep}
Soit $\Gamma$ un sous-groupe discret de $G$.
\begin{enumerate}
\item La vari\'et\'e complexe $\Gamma\backslash S$ est holomorphiquement
s\'eparable.
\item Si $G/\Gamma$ est compact, alors les fonctions holomorphes born\'ees
s\'eparent les points de $\Gamma\backslash S$, i.\,e.\ $\Gamma\backslash S$ est
hyperbolique au sens de Caratheodory.
\item Si $G={\rm{SL}}(2,\mbb{R})$ et si $\Gamma$ est contenu dans
${\rm{SL}}(2,\mbb{Z})$, alors $\Gamma\backslash S$ est une vari\'et\'e de Stein.
\end{enumerate}
\end{thm}

La preuve de ce th\'eor\`eme utilise la th\'eorie des repr\'esentations de $S$.

Enfin, nous notons l'observation suivante qui se trouve sans d\'emonstration
dans~\cite{BrRu}. Pour l'aise du lecteur nous en donnons une d\'emonstration.
\`A l'egard de la th\'eorie des ensembles r\'eel-alg\'ebriques et
semi-alg\'ebriques nous renvoyons le lecteur \`a~\cite{BCR}.

\begin{lem}\label{Lem:semialgebraic}
Le semi-groupe $S$ est un sous-ensemble semi-alg\'ebrique de $G^\mbb{C}$.
\end{lem}

\begin{proof}
Soit $A\colon G^\mbb{C}\times Z\to Z\times Z$ l'application $(g,z)\mapsto(g
\cdot z,z)$ et soit $p_{G^\mbb{C}}\colon G^{\mbb{C}}\times Z\to G^\mbb{C}$ la
projection sur la premi\`ere composante. Il n'est pas difficile de montrer que
l'on a:
\begin{equation*}
G^\mbb{C}\setminus S=p_{G^\mbb{C}}\bigl(A^{-1}(Z\setminus
X\times\overline{X})\bigr).
\end{equation*}

De plus il est bien connu que $G^\mbb{C}$ porte une unique structure de groupe
lin\'eaire-alg\'ebrique et que $G^\mbb{C}$ agit de fa\c con alg\'ebrique
sur la vari\'et\'e projective $Z$. Puisque $G$ est r\'eel-alg\'ebrique dans
$G^\mbb{C}$, le th\'eor\`eme de Tarski-Seidenberg implique que l'orbite
$X=G\cdot eQ^-$ est semi-alg\'ebrique. Par cons\'equent, son compl\'ementaire
$Z\setminus X$ et son adh\'erence $\overline{X}$ sont \'egalement
semi-alg\'ebriques, donc $A^{-1}(Z\setminus X\times\overline{X})$ est
semi-alg\'ebrique. Enfin, le th\'eor\`eme de Tarski-Seidenberg implique que
$G^\mbb{C}\setminus S$ et par cons\'equent $S$ sont semi-alg\'ebriques.
\end{proof}

\section{Une condition suffisante pour que $\Gamma\backslash S$ soit de Stein}

Nous continuons la notation etablie dans la section pr\'ec\'edente. La
complexification de l'espace symmetrique Hermitien $X=G/K$ est par d\'efinition
la vari\'et\'e complexe homog\`ene $X^\mbb{C}:=G^\mbb{C}/K^\mbb{C}$. Puisque
$K^\mbb{C}$ est ferm\'e dans $Q^-\cong K^\mbb{C}\ltimes P^-$ et puisque $Q^-$
agit effectivement sur $P^-$, l'application naturelle $p\colon G^\mbb{C}/
K^\mbb{C}\to G^\mbb{C}/Q^-$ est un fibr\'e holomorphe dont la fibre est donn\'ee
par $P^-\cong\mbb{C}^{\dim X}$ et dont le groupe structural est le groupe
complexe connexe $Q^-$.

Il est connu que l'application $gK^\mbb{C}\mapsto(gQ^-,gQ^+)$ est un
plongement ouvert de $X^\mbb{C}$ dans $G^\mbb{C}/Q^-\times G^\mbb{C}/Q^+$. Il
s'ensuit que nous obtenons que le diagramme
\begin{equation*}
\xymatrix{
G^\mbb{C}/K^\mbb{C}\ar[r]\ar[d] & G^\mbb{C}/Q^-\times G^\mbb{C}/Q^+\ar[d]\\
G^\mbb{C}/Q^-\ar[r] & G^\mbb{C}/Q^-
}
\end{equation*}
est commutatif.

Soit $\pi\colon G^\mbb{C}\to G^\mbb{C}/K^\mbb{C}$ le fibr\'e principal de groupe
structural $K^\mbb{C}$. Nous commen\c cons par montrer que l'image $\pi(S)$ est
une vari\'et\'e de Stein. Un th\'eor\`eme de Matsushima et Morimoto assure que
$\pi(S)$ est de Stein si et seulment si $SK^\mbb{C}$ est de Stein.

\begin{prop}\label{Prop:SteinGlob}
Le domaine $SK^\mbb{C}\subset G^\mbb{C}$ est de Stein. Par cons\'equent,
$\pi(S)\subset G^\mbb{C} /K^\mbb{C}$ est de Stein.
\end{prop}

\begin{proof}
Puisque $S$ est invariant par rapport \`a la multiplication \`a droit de $K$,
on obtient une action holomorphe locale par $K^\mbb{C}$ sur $S$. Par~\cite{He}
il y a une vari\'et\'e de Stein $S^*$ munie d'une action holomorphe de
$K^\mbb{C}$ et un plongement ouvert holomorphe $K$--equivariant $\varphi\colon
S\to S^*$ tel que $S^*=K^\mbb{C}\cdot\varphi(S)$. De plus, $S^*$ est universel
ce qui implique l'existence d'une application holomorphe
$K^\mbb{C}$--equivariante $\Phi\colon S^*\to SK^\mbb{C}$ tel que
$\Phi\circ\varphi$ co\"\i ncide avec l'inclusion $S \hookrightarrow
SK^\mbb{C}$. Il s'en suit que $\Phi$ est surjectif et localement biholomorphe.
Donc il suffit de montrer que $\Phi$ est une application finie. Dans notre
situation c'est le cas si et seulement si les fibres de $\Phi$ sont finies.

Pour $z_0\in S$ nous d\'efinisson l'ensemble
\begin{equation*}
K^\mbb{C}[z_0]:=\bigl\{k\in K^\mbb{C};\ z_0k^{-1}\in S\bigr\}.
\end{equation*}
Or, $S$ est semi-alg\'ebrique par Lemme~\ref{Lem:semialgebraic} et $K^\mbb{C}$
est un groupe lin\'eaire alg\'ebrique, donc $K^\mbb{C}[z_0]$ est
semi-alg\'ebrique. En particulier, le nombre de composantes connexes de
$K^\mbb{C}[z_0]$ est fini. Nous allons d\'emontrer que $\#\Phi^{-1}(z_0)$ est
born\'e par le nombre de composantes connexes de $K^\mbb{C}[z_0]$.

Soit $z\in\Phi^{-1}(z_0)$. Alors il y a un \'el\'ement $k\in K^\mbb{C}$ tel que 
$k\cdot z\in\varphi(S)$. Par cons\'equent, nous avons $\Phi(k\cdot z)=
k\cdot\Phi(z)=z_0k^{-1}\in S$, i.\,e.\ $k\in K^\mbb{C}[z_0]$. Comme $\Phi
|_{\varphi(S)}\colon\varphi(S)\to S$ est biholomorphe, nous concluons $k\cdot
z=\varphi(z_0k^{-1})$, ce qui est \'equivalent \`a
$z=k^{-1}\cdot\varphi(z_0k^{-1})$. Comme r\'eciproquement tout \'el\'ement $z\in
S^*$ de la forme $z=k^{-1}\cdot\varphi(z_0k^{-1})$ avec $k\in K^\mbb{C}[z_0]$
est contenu dans $\Phi^{-1}(z_0)$, nous obtenons une application surjective
\begin{equation}\label{Eqn:Map}
K^\mbb{C}[z_0]\to\Phi^{-1}(z_0),\quad k\mapsto k^{-1}\cdot\varphi(z_0k^{-1}).
\end{equation}
Si $k\in K^\mbb{C}[z_0]$ est contenu dans la composante connexe de l'\'el\'ement
neutre, nous avons $\varphi(z_0k^{-1})=k\cdot\varphi(z_0)$. Cette
observation implique que l'application~\eqref{Eqn:Map} est localement constante,
ce qui montre la proposition.
\end{proof}

\begin{rem}
La vari\'et\'e $S^*$ est la globalisation universelle de la $K^\mbb{C}$--action
locale sur $S$ au sens de~\cite{Pa}.
\end{rem}

La proposition~\ref{Prop:SteinGlob} nous permet de donner la condition
equivalente suivante pour que $\Gamma\backslash S$ soit de Stein.

\begin{prop}\label{Prop:KCQuot}
Soit $\Gamma\subset G$ un sous-groupe discret. Alors $\Gamma\backslash S$ est
de Stein si et seulement si $\Gamma\backslash\pi(S)$ est de Stein.
\end{prop}

\begin{proof}
Si $\Gamma\backslash S$ est de Stein, alors la globalisation universelle $\Gamma
\backslash S^*$ de la $K^\mbb{C}$--action locale sur $\Gamma\backslash S$ est
aussi de Stein. Comme $K^\mbb{C}[\Gamma g]=K^\mbb{C}[g]$ pour tout $g\in S$, le
m\^eme argument qu'au-dessus implique que $\Gamma\backslash SK^\mbb{C}$ est
de Stein, ce qui montre que $\Gamma\backslash\pi(S)$ est de Stein.

R\'eciproquement, si $\Gamma\backslash\pi(S)$ est de Stein, alors
$\Gamma\backslash SK^\mbb{C}$ est aussi de Stein. Puisque $\Gamma\backslash S$
est plong\'e comme ouvert dans $\Gamma\backslash SK^\mbb{C}$ et localement de
Stein, l'affirmation d\'ecoule d'un argument \`a la Docquier-Grauert (comparer
Proposition~4.7 dans~\cite{Mie4}).
\end{proof}

Soit $p\colon G^\mbb{C}/K^\mbb{C}\to G^\mbb{C}/Q^-$ le fibr\'e holomorphe et
consid\'erons $X\subset G^\mbb{C}/Q^-$.

\begin{lem}
On a $\pi(S)\subset p^{-1}(X)$.
\end{lem}

\begin{proof}
Soit $g\in S$. Comme les applications $\pi$ and $p$ sont
$G^\mbb{C}$--equivariantes, nous en d\'eduisons que $p\circ\pi(g)=g\cdot
p\circ\pi(e)=g\cdot eQ^-$. \`A cause de $eQ^-\in X$, la d\'efinition de $S$
implique que $gQ^-\in X$.
\end{proof}

La restriction $p\colon p^{-1}(X)\to X$ est un fibr\'e holomorphe de fibre $P^-$
et de groupe structural complexe connexe $Q^-$. De plus, le groupe $G$ agit
par des automorphismes de fibr\'e sur $p^{-1}(X)$.

Le th\'eor\`eme suivant donne une condition suffisante pour que
$\Gamma\backslash S$ soit de Stein.

\begin{thm}
Soit $\Gamma\subset G$ un sous-group discret qui agit librement sur $X$. Si
$X/\Gamma$ est de Stein, alors $\Gamma\backslash S$ est \'egalement de Stein.
\end{thm}

\begin{proof}
Comme le groupe $\Gamma$ agit par automorphismes de fibr\'e, nous obtenons
le fibr\'e holomorphe $\Gamma\backslash p^{-1}(X)\to X/\Gamma$ (comparer
Proposition~6.1 dans~\cite{Mie4}). Comme $X/\Gamma$ est de Stein, la fibre
$P^-$ est de Stein, et le groupe structural de ce fibr\'e est complexe
connexe, nous en d\'eduisons que $\Gamma\backslash p^{-1}(X)$ est de Stein.
Puisque $\pi(S)$ est de Stein et $\Gamma\backslash\pi(S)$ est un ouvert dans
$\Gamma\backslash p^{-1}(X)$, nous concluons que $\Gamma\backslash\pi(S)$ est
de Stein et donc que $\Gamma\backslash S$ est de Stein.
\end{proof}

\begin{ex}
\begin{enumerate}
\item Soit $G={\rm{SL}}(2,\mbb{R})$ et soit $\Gamma\subset G$ un sous-groupe
discret qui agit librement sur le demi-plan sup\'erieur $\mbb{H}$ tel que la
surface de Riemann $R:=\mbb{H}/\Gamma$ est non-compacte. Dans ce cas $R$
est de Stein, donc le quotient $\Gamma\backslash S$ est de Stein aussi. Cette
observation donne une nouvelle d\'emonstration du
Th\'eor\`eme~\ref{Thm:HolSep}~(3).
\item Si $X=G/K$ est un espace symmetrique hermitien arbitraire de type
non-compact et si $\Gamma\subset G$ et cyclique, alors $X/\Gamma$ et donc
$\Gamma\backslash S$ sont de Stein par~\cite{Mie4}.
\end{enumerate}
\end{ex}

\section{Un contre-exemple}

Dans cette section nous consid\'erons l'exemple explicit du disque $X=\Delta=\{
z\in\mbb{C};\ \abs{z}<1\}$. Le dual compact de $X$ est donn\'e par $Z=\mbb{P}_1$
et nous avons
\begin{equation*}
\Delta\cong\bigl\{[z_0:z_1]\in\mbb{P}_1;\ \abs{z_0}^2-\abs{z_1}^2<0\bigr\}\cong
G/K
\end{equation*}
pour $G={\rm{SU}}(1,1)$ et $K=\left\{\left(\begin{smallmatrix}
e^{it}&0\\0&e^{-it}\end{smallmatrix}\right);\ t\in\mbb{R}\right\}\cong S^1$.

On v\'erifie que $P^-=\left(\begin{smallmatrix}1&0\\\mbb{C}&1
\end{smallmatrix}\right)$ et $P^+=\left(\begin{smallmatrix}
1&\mbb{C}\\0&1\end{smallmatrix}\right)$. Donc, nous pouvons identifier
$G^\mbb{C}/Q^-\times G^\mbb{C}/Q^+$ avec $\mbb{P}_1\times\mbb{P}_1$ o\`u le
point $(eQ^-,eQ^+)$ correspond \`a $\bigl([0:1],[1:0]\bigr)$. De plus, nous
identifions la complexification $X^\mbb{C}=G^\mbb{C}/K^\mbb{C}$ avec la
$G^\mbb{C}$--orbite ouverte de $(eQ^-,eQ^+)$, i.\,e.\ avec $(\mbb{P}_1\times
\mbb{P}_1)\setminus D$ o\`u nous notons par $D$ la diagonale dans
$\mbb{P}_1\times\mbb{P}_1$. Par cons\'equent, $p^{-1}(\Delta)=(\Delta\times
\mbb{P}_1)\setminus D$ contient les domaines $G$--invariants
$(\Delta\times\Delta)\setminus D$ et
$\Delta\times(\mbb{P}_1\setminus\overline{\Delta})$ qui sont de Stein tous les
deux. Le groupe $G$ agit librement et transitivement sur leur fronti\`ere
$F:=\Delta\times S^1$ qui est alors une hypersurface Levi-plate dans
$p^{-1}(\Delta)$.

\begin{lem}
Le domaine $\pi(S)$ contient l'hypersurface Levi-plate $F$.
\end{lem}

\begin{proof}
On v\'erifie directement que pour $t>0$ l'\'el\'ement $\left(\begin{smallmatrix}
e^{-t}&0\\0&e^t\end{smallmatrix}\right)$ est contenu dans $S\cap K^\mbb{C}$.
Comme la matrice $\left(\begin{smallmatrix} 1+ix&-ix\\ix&1-ix
\end{smallmatrix}\right)$ se trouve dans $G={\rm{SU}}(1,1)$ pour tout
$x\in\mbb{R}$, leur produit est contenu dans $S$. L'orbite par ce produit du
point $\bigl([0:1],[1:0]\bigr)$ est donne\'e par
\begin{equation*}
\bigl([-ie^{-t}x:e^t(1-ix)],[e^{-t}(1+ix):ie^tx]\bigr).
\end{equation*}
Un calcul simple montre que pour le choix $x=1$ et $t=\tfrac{1}{4}\log(2)$ cet
\'el\'ement est contenu dans $F$. En utilisant la multiplication \`a gauche par
$G$ l'affirmation est demontr\'ee.
\end{proof}

\begin{prop}
Si $G/\Gamma$ est compact, alors le quotient $\Gamma\backslash S$ n'est pas de
Stein.
\end{prop}

\begin{proof}
Par la proposition~\ref{Prop:KCQuot} le quotient $\Gamma\backslash S$ est de
Stein si et seulement si $\Gamma\backslash\pi(S)$ est de Stein. Nous avons vu
que $\pi(S)$ contient une hypersurface $G$--invariante Levi-plate $F$ ce qui
implique que $\Gamma\backslash\pi(S)$ contient la hypersurface compacte
Levi-plate $\Gamma\backslash F$. Soit $f$ une fonction holomorphe arbitraire sur
$\Gamma\backslash\pi(S)$ et soit $z\in\Gamma\backslash F$ un point dans lequel
la norme absolue de la restriction de $f$ sur $\Gamma\backslash F$ prend son
maximum. Comme $\Gamma\backslash F$ est Levi-plate, il y a une feuille
holomorphe qui passe par $z$. Le principe du maximum implique que $f$ doit
\^etre constante sur cette feuille. Nous en d\'eduisons que les points de cette
feuille ne peuvent pas \^etre s\'epar\'es par les fonctions holomorphes, i.\,e.\
que $\Gamma\backslash\pi(S)$ n'est pas holomorphiquement s\'eparable, et en
particulier, pas de Stein.
\end{proof}

\end{document}